\documentclass[12pt]{article}
\usepackage{times,amsmath,amsthm,amssymb,graphicx,xspace,epsfig,xcolor}
\usepackage[plain]{fullpage}
\usepackage{thm-restate}
\usepackage{hyperref}

\usepackage{tikz}

\tikzstyle{vertex}=[circle,draw, top color=white, 
	    bottom color=white, minimum size=20pt, scale=1, inner sep=0.5pt]
\tikzstyle{edge}=[thick]
\tikzstyle{arc}=[->, thick]

\usepackage{color}
\usepackage{authblk,enumitem}

\newtheorem{theorem}{Theorem}[section]

\newtheorem{proposition}[theorem]{Proposition}
\newtheorem{lemma}[theorem]{Lemma}
\newtheorem{claim}{Claim}[theorem]

\theoremstyle{definition}
\newtheorem{remark}[theorem]{Remark}

\newtheorem{problem}[theorem]{Problem}
\newtheorem{conjecture}[theorem]{Conjecture}

\newcommand{\rev}[1]{\textcolor{black}{#1}}

\definecolor{dark-green}{rgb}{0.2, 0.5, 0.2}

\newcommand{\ra}{\rightarrow}
\newcommand{\Ra}{\Rightarrow}

\newenvironment{subproof}{\par\noindent {\it Proof of claim}.\ }{\hfill$\lozenge$\par\vspace{11pt}}

\DeclareMathOperator{\rk}{rk}
\DeclareMathOperator{\Inv}{Inv}

\DeclareMathOperator{\inv}{inv}

\DeclareMathOperator{\minv}{minv}

\DeclareMathOperator{\bigO}{\mathcal{O}}

\title{Problems, proofs, and disproofs on the inversion number}

\author[1,2]{Guillaume Aubian}
\author[3]{Fr\'ed\'eric Havet}
\author[4]{Florian H\"orsch}
\author[5]{Felix Klingelhoefer}
\author[3]{Nicolas Nisse}
\author[2,3]{Cl\'ement Rambaud}
\author[2]{Quentin Vermande}

\affil[1]{Universit\'e de Paris, CNRS, IRIF, Paris, France}
\affil[2]{DIENS, \'Ecole normale sup\'erieure, CNRS, PSL University, Paris, France}
\affil[3]{Universit\'e C\^ote d'Azur, CNRS, Inria, I3S, Sophia Antipolis, France}
\affil[4]{Technische Universit\"at Ilmenau, Ilmenau, Germany}
\affil[5]{Universit\'e Grenoble Alpes, Grenoble, France}

\date{}

\begin{document}
\maketitle

\begin{abstract}
The {\it inversion} of a set $X$ of vertices in a digraph $D$ consists in reversing the direction of all arcs of $D\langle X\rangle$. The {\it inversion number} of an oriented graph $D$, denoted by $\inv(D)$,  is the minimum number of inversions needed to transform $D$ into an acyclic oriented graph. In this paper, we study a number of problems involving the inversion number of oriented graphs.
Firstly, we give bounds on $\inv(n)$, the maximum of the inversion numbers of the oriented graphs of order $n$.
We show $n - \bigO(\sqrt{n\log n}) \ \leq \ \inv(n) \ \leq \  n - \lceil \log (n+1) \rceil$. 
Secondly, we disprove a conjecture of Bang-Jensen et al.~\cite{BCH} asserting that, for every pair of oriented graphs $L$ and $R$, we have $\inv(L\Ra R) =\inv(L) + \inv(R)$, where $L\Ra R$ is the oriented graph obtained from the disjoint union of $L$ and $R$ by adding all arcs from $L$ to $R$. Finally, we investigate whether, for all pairs of positive integers $k_1,k_2$, there exists an integer $f(k_1,k_2)$ such that if $D$ is an oriented graph with $\inv(D) \geq f(k_1,k_2)$ then there is a partition $(V_1, V_2)$ of $V(D)$ such that $\inv(D\langle V_i\rangle) \geq k_i$ for $i=1,2$. We show that $f(1,k)$ exists and $f(1,k)\leq k+10$ for all positive integers $k$. Further, we show that $f(k_1,k_2)$ exists for all pairs of positive integers $k_1,k_2$ when the oriented graphs in consideration are restricted to be tournaments.
\medskip

    \noindent{}{\bf Keywords:}  inversion; tournament; oriented graph.
  \end{abstract}

\section{Introduction}

Notation not given below is consistent with \cite{bang2009}. We denote by $[k]$ the set $\{1,2, \dots, k\}$. The notation $\log$ refers to the logarithm to the base 2.

Let $D$ be an oriented graph.
The {\bf inversion} of a set $X$ of vertices of $D$ consists in reversing the direction of all arcs of $D\langle X\rangle$.
We say that we {\bf invert} $X$ in $D$. The resulting oriented graph is denoted by $\Inv(D;X)$.  
If $(X_i)_{i\in I}$  is a family of subsets of $V(D)$, then $\Inv(D; (X_i)_{i\in I})$ is the oriented graph obtained after inverting the
$X_i$ one after another. Observe that this is independent of the order in which we invert the $X_i$~: $\Inv(D; (X_i)_{i\in I})$ is obtained from $D$ by reversing the arcs such that an odd number of the $X_i$ contain its two end-vertices.
A {\bf decycling family} of an oriented graph $D$ is a family of subsets $(X_i)_{i\in I}$ of subsets of $V(D)$ such that
 $\Inv(D; (X_i)_{i\in I})$ is acyclic.
 The {\bf inversion number} of an oriented graph $D$, denoted by $\inv(D)$,  is the minimum number of inversions needed to transform $D$ into an acyclic oriented graph, that is, the minimum cardinality of a decycling family.

The inversion number of oriented graphs was first introduced by Belkhechine et al. in \cite{BBBP10}.
They studied oriented graphs with large inversion number.
For any positive integer $n$, let $\inv(n)=\max\{\inv(D) \mid D~\mbox{oriented graph of order}~n\}$. Since the inversion number is non-increasing with respect to edge deletion, we have
$\inv(n)=\max\{\inv(T) \mid T~\mbox{tournament of order}~n\}$. 
In the following, we give some basic bounds on $\inv(n)$.

\begin{remark} \label{rem:inv}
 $\inv(n) \leq \inv(n-1)+1$ for all positive integer $n$.
\end{remark}
\begin{proof}
Let $T$ be a tournament of order $n$.
Pick a vertex $x$ of $T$. Observe that $x$ is a sink in $T'=\Inv(T; N_T^+[x]\cup x)$. 
This yields $\inv(T') =\inv(T'-x) \leq \inv(n-1)$. Hence $\inv(T) \leq \inv(n-1)+1$.
\end{proof}

Every oriented graph on at most two vertices is acyclic, so $\inv(1)=\inv(2)=0$.
Every tournament of order at most $4$ has a cycle arc-transversal of size at most $1$, so $\inv(3)=\inv(4)=1$.
As observed by Belkhechine et al.~\cite{BBBP}, every tournament of order at most $6$ has inversion number at most $2$.  This observation and Remark \ref{rem:inv} yield

\begin{equation}\label{up-easy}
\inv(n) \leq n-4  ~~~~\mbox{for all}~n\geq 6.
\end{equation}

Moreover, Belkhechine et al.~\cite{BBBP10} observed that since there are $n!$ labelled transitive tournaments of order $n$, the number of labelled tournaments of order $n$ with inversion number less than $p$ is at most $n! 2^{n(p-1)}$, while there are $2^{\frac{n(n-1)}{2}}$ labelled tournaments of order $n$. It follows that for all positive integers $n$ and $p$ such that $2^{\frac{n(n-1)}{2}} >  n! 2^{n(p-1)}$, there is a tournament $T$ of order $n$ such that $\inv(T) \geq p$. Hence, for every positive integer $n$, we have 

\begin{equation}\label{down-easy}
\inv(n) \geq \frac{n-1}{2} - \log n.  ~~~~
\end{equation}

However, Belkhechine et al.~\cite{BBBP} conjectured that Equation~\eqref{down-easy} is not tight. 

 \begin{conjecture}[Belkhechine et al.~\cite{BBBP}]\label{conj:invn}
 $\inv(n) \geq \lfloor\frac{n-1}{2}\rfloor$.
 \end{conjecture}

 Their conjecture was motivated by their conviction that some explicit tournament of order $n$ has inversion number at least  $\lfloor\frac{n-1}{2}\rfloor$. 
 Let $Q_n$ be the tournament obtained from the transitive tournament by reversing the arcs of its unique directed hamiltonian path $(v_1,v_2, \ldots{},v_n)$.

 \begin{conjecture}[Belkhechine et al.~\cite{BBBP}]\label{conj:Q}
 $\inv(Q_n)= \lfloor\frac{n-1}{2}\rfloor$.
 \end{conjecture}

 The first main contribution of this article is an improvement of both the upper and the lower bound on $\inv(n)$. With a rather simple argument we improve the upper bound.
 
\begin{restatable}{theorem}{invupper}\label{thm:inv-upper}
For every positive integer $n$, we have $\inv(n) \leq n - \lceil \log (n+1) \rceil$.
\end{restatable}

Next, we improve the lower bound using some more involved probabilistic techniques. More precisely, considering a random tournament of order $n$, we show the following result.

\begin{restatable}{theorem}{random}\label{theorem:random_tournament}
$\inv(n) \geq n - 2\sqrt{n\log n}$ for $n$ sufficiently large.
\end{restatable}

Observe that for $n$ sufficiently large, Theorem \ref{theorem:random_tournament} confirms Conjecture \ref{conj:invn} and even establishes a much stronger bound. The proofs of Theorems \ref{thm:inv-upper} and \ref{theorem:random_tournament} can be found in Section \ref{sec:invn}.


\bigskip

A {\bf cycle transversal} (resp.  {\bf cycle arc-transversal}) in a digraph $D$ is a set of vertices (resp. arcs) whose deletion results in an acyclic digraph. We wish to remark that cycle transversals and cycle arc-transversals are sometimes also referred to in the literature as feedback vertex sets and feedback arcs sets, respectively.
The {\bf cycle transversal number} (resp. {\bf  cycle arc-transversal number}) is the minimum size of a cycle transversal (resp. cycle arc-transversal) of $D$ and is denoted by $\tau(D)$ (resp. $\tau'(D)$).  
As observed in~\cite{BBBP}, the inversion number is bounded by the cycle arc-transversal number and twice the cycle transversal number.

\begin{theorem}[Bang-Jensen et al.~\cite{BCH}]\label{thm:bound-fvs}
$\inv(D) \leq  \tau' (D)$ and  $\inv(D) \leq 2\tau(D)$ for every oriented graph $D$.
 \end{theorem}

A natural question is to ask whether these bounds are tight or not.

We denote by $\vec{C_3}$ the directed cycle of length $3$ and by $TT_n$ the transitive tournament of order $n$. The vertices of $TT_n$ are $v_1, \dots, v_n$ and its arcs $\{v_iv_j \mid i < j\}$.
The {\bf lexicographic product}  of a digraph $D$ by a digraph $H$ is the digraph $D[H]$ with vertex set $V(D)\times V(H)$
and arc set
$A(D[H])  =  \{(a,x)(b,y) \mid ab \in A(D), \mbox{\ or\ } a=b \mbox{\ and\ } xy\in A(H)\}$.
It can be seen as blowing up each vertex of $D$ by a copy of $H$.
Using boolean dimension,  Pouzet et al.~\cite{PST} proved  the following.
\begin{theorem}[Pouzet et al.~\cite{PST}]\label{thm:TT[C3]}
 $\inv(TT_n[\vec{C_3}]) = n$. 
\end{theorem}
Since $\tau'(TT_n[\vec{C_3}])=n$, this shows that the inequality $\inv(D) \leq  \tau' (D)$  of Theorem~\ref{thm:bound-fvs} is tight.

Pouzet asked for an elementary proof of Theorem~\ref{thm:TT[C3]}.
Let $L$ and $R$ be two oriented graphs. The {\bf dijoin} from $L$ to $R$, denoted by $L\Ra R$, is the oriented graph obtained from the disjoint union of $L$ and $R$ by adding all arcs from $L$ to $R$.
Observe that $TT_n[\vec{C_3}]=\vec{C_3}\Ra TT_{n-1}[\vec{C_3}]$. Hence one way to prove Theorem~\ref{thm:TT[C3]} would be to prove that holds $\inv(\vec{C_3}\Ra T) = \inv(T)+1$ for every tournament $T$. First inverting $\inv(L)$ subsets of $V(L)$ to make $L$ acyclic and then inverting $\inv(R)$ subsets of $V(R)$ to make $R$ acyclic, makes $L\Ra R$ acyclic. Therefore $\inv(L\Ra R) \leq \inv(L) + \inv(R)$.
Bang-Jensen et al. \cite{BCH} conjectured that equality always holds.

\begin{conjecture}[Bang-Jensen et al. \cite{BCH}]\label{conj:dijoin}
For any two oriented graphs, $L$ and $R$, $\inv(L\Ra R) =\inv(L) + \inv(R)$.
\end{conjecture}

Bang-Jensen et al.~\cite{BCH} showed that the inequality  $\inv(D) \leq 2\tau(D)$ of Theorem~\ref{thm:bound-fvs} is tight for $\tau(D)\in \{1,2\}$ and conjectured the following.

\begin{conjecture}[Bang-Jensen et al.~\cite{BCH}] \label{conj:tight}
 For every positive integer $n$, there exists an oriented graph $D$ such that $\tau (D)=n$ and  $\inv(D) =2n$.
 \end{conjecture}

 This conjecture would also be implied by Conjecture~\ref{conj:dijoin}.
 
Moreover, Bang-Jensen et al.~\cite{BCH} proved that deciding whether a given digraph $D$ has inversion number at most $1$ is NP-complete. Together with Conjecture~\ref{conj:dijoin}, this would imply the following conjecture posed in ~\cite{BCH}.

 \begin{conjecture}[Bang-Jensen et al.~\cite{BCH}]\label{conj:NP}
 Deciding whether a given digraph $D$ has inversion number at most $k$ is NP-complete for any fixed positive integer $k$. \end{conjecture}

 \medskip

Bang-Jensen et al.~\cite{BCH} proved that Conjecture~\ref{conj:dijoin} holds when $\inv(L)\leq 1$ and $\inv(R)\leq 2$, and when $\inv(L)=\inv(R)=2$ and both $L$ and $R$ are strongly connected.
Unfortunately, in Section~\ref{sec:dijoin}, we disprove Conjecture~\ref{conj:dijoin}. More precisely, we show the following result:
\begin{restatable}{theorem}{disproof}\label{cor:disproof}
For every odd integer $k\geq 3$, there is a tournament $T_k$ with  $\inv(T_k)=k$ such that
$\inv(T_k \Ra R) \leq k+\inv(R) -1$ for every oriented graph $R$ with $\inv(R)\geq 1$.
\end{restatable}

\bigskip
A celebrated conjecture of Alon~\cite{Alon96} states that for any two positive integers $k_1,k_2$, there is 
a minimum integer $f(k_1,k_2)$ such that, for any digraph $D$ with minimum out-degree $f(k_1,k_2)$, there is a partition $(V_1, V_2)$ of $V(D)$ such that the minimum out-degree of the subdigraph induced by $V_i$ is at least $k_i$ for each $i=1,2$.
It is natural to ask whether an Alon-type result holds for the inversion number. 

\begin{conjecture}\label{pb:split}
For every two positive integers $k_1,k_2$, there exists an integer $f(k_1,k_2)$ such that every oriented graph $D$ with $\inv(D) \geq f(k_1,k_2)$ admits a partition $(V_1, V_2)$ of $V(D)$ such that $\inv(D\langle V_i\rangle) \geq k_i$ for $i=1,2$.
\end{conjecture}

An oriented graph is {\bf intercyclic} if it does not have two vertex-disjoint directed cycles.
Bang-Jensen et al.~\cite{BCH} proved that intercyclic oriented graphs have inversion number at most $4$. This implies that $f(1,1)$ exists and $f(1,1) \leq 5$.

In an attempt to approach Conjecture~\ref{pb:split}, we give two partial results. Firstly, we show that an analogous statement holds when we restrict the oriented graphs in consideration to be tournaments.

\begin{restatable}{theorem}{ftour}\label{cor:f-tournament}
For every two positive integers $k_1,k_2$, there exists an integer $f_T(k_1,k_2)$ such that every tournament $T$ with $\inv(T) \geq f_T(k_1,k_2)$ admits a partition $(V_1, V_2)$ of $V(T)$ such that $\inv(T\langle V_i\rangle) \geq k_i$ for $i=1,2$.
\end{restatable}

Secondly, we show that the statement holds when restricting one of the two parameters to be 1. More concretely, we show the following result.

\begin{theorem}\label{thm:1split}
For every positive integer $k$ and for every oriented graph $D$ with $\inv(D)\geq  k+10$, there is a partition $(V_1, V_2)$ of $V(D)$ satisfying $\inv(D\langle V_1\rangle) \geq k$ and $\inv(D\langle V_2\rangle) \geq 1$.
\end{theorem} The proofs of Theorems~\ref{cor:f-tournament} and \ref{thm:1split} can be found in Section \ref{sec:splitting}.

\medskip

We conclude this work in Section~\ref{sec:further}, by proposing directions for further research on the topic.

First, since the lower bound on $\inv(n)$ of Theorem~\ref{theorem:random_tournament} is obtained by random tournaments, we explore possible explicit constructions of tournaments with large inversion number.
Theorem~\ref{thm:TT[C3]} yields tournaments of order $n$ with inversion number $\lfloor n/3\rfloor$. 
Explicit constructions of tournaments with larger inversion number are not known.
Proving (or approaching) Conjecture~\ref{conj:Q} would give one.
Let us describe another construction.
Let $D$ be an oriented graph.
We denote by $\triangle(D)$ the oriented graph obtained from $D$ by adding two new vertices $u,v$ and all the arcs from $V(D)$ to $u$, $uv$, and all the arcs from $v$ to $V(D)$. We believe that iterating this construction yields tournaments with large inversion number.

\begin{restatable}{conjecture}{conjtriangle}\label{conj:triangle}
 $\inv(\triangle^{k}(TT_1)) = \inv(\triangle^{k}(TT_2)) = k+1$.
\end{restatable} 

In Subsection~\ref{subsec:explicit}, we consider this conjecture. In particular, in Theorem~\ref{thm:triangle-linear}, we prove that
$\inv(\triangle^{k-1}(TT_2)) \geq k/2-1$. 

Finally, in Subsection~\ref{subsec:disproof}, we
present some problems raised by the disproof of Conjecture~\ref{conj:dijoin}.

\paragraph{Note added}
Almost simultaneously with the release of this paper, Alon, Powierski, Savery, Scott, and Wilmer announced independent work~\cite{APSSW}
on some of the problems we address here. Specifically, they show upper and lower bounds on
$\inv(n)$ of forms similar to those of Theorems~\ref{thm:inv-upper}  and \ref{theorem:random_tournament}. 
They also disprove Conjecture~\ref{conj:dijoin}, but only provide one counterexample. In addition, they prove Conjectures~\ref{conj:tight} and~\ref{conj:NP}.

\section{Maximum inversion number of an oriented graph of order $n$} \label{sec:invn}

This section is dedicated to improving the bounds on $\inv(n)$. More particularly, we prove Theorems~\ref{thm:inv-upper} and \ref{theorem:random_tournament}.

We first give the proof of Theorem \ref{thm:inv-upper} which we recall.

\invupper*

\begin{proof}
As observed above, it suffices to prove the statement for tournaments. We proceed
by induction on $n$, the result holding trivially when $n = 1$.    

Let $T$ be a tournament on $n$ vertices for some $n \geq 2$.
Now consider a vertex $u$. Without loss of generality, we may assume that $d_T^-(u) \leq \lfloor\frac{n-1}{2}\rfloor$.
Let $v_1, \dots , v_{d_T^-(u)}$ be an arbitrary ordering of $N_T^-(u)$.
For every $i \in [d_T^-(u)]$, let $X_i = \{v_i\} \cup N^-_{T_{i-1}}(v_i) \setminus \{v_j \mid j < i\}$ where $T_{i-1}=\Inv(T; (X_j)_{j\in[i-1]})$.
Observe that after inverting $X_i$, $v_i$ dominates all the vertices of $\{v_j \mid j > i\} \cup \{u\} \cup N^+_T(u)$,
and $u$ dominates $N^+_T(u)$.
Hence it remains to apply the induction hypothesis to make the subtournament induced by $N_T^+(u)$ acyclic.
In total, we used at most $\lfloor\frac{n-1}{2}\rfloor+\lceil\frac{n-1}{2}\rceil-\log(\lceil\frac{n-1}{2}\rceil+1) = n - \log(2\lceil\frac{n-1}{2}\rceil+2)
\leq n - \log(n+1)$ inversion. This proves the theorem.
\end{proof}


\medskip

The proof of Theorem \ref{theorem:random_tournament} is based on probabilistic methods and heavily relies on properties of a certain random matrix over the unique field $\mathbb{F}_2$ on two elements 0 and 1. All the matrices considered in this paper are over this field, and for a matrix $M$, we denote by $\rk(M)$ its rank over $\mathbb{F}_2$. A significant part of the technicalities of the proof is included in the following lemma.

\begin{lemma}\label{lem:random_matrix}
Let $n \geq 1$ be a positive integer and let $x_1, \dots, x_n \in \mathbb{F}_2$.
Let $M$ be a matrix in $\mathbb{F}_2^{n \times n}$ chosen uniformly
at random among all symmetric matrices whose diagonal entries are $x_1, \dots, x_n$.
Then for every $r \in \{0,\ldots,n\}$, we have
\[
\Pr[\rk(M) \leq n-r] \leq  2^{-\frac{1}{2}(r^2-4r)}
\]
\end{lemma}
\begin{proof}
For $i \in \{0,\ldots,n\}$, let $M_i$ be the submatrix of $M$ which is restricted to the first $i$ rows and columns. Observe that $M_n=M$. Further, we define random variables $X_i,i=0,\ldots,n$ where $X_i=i-\rk(M_i)$.
\begin{claim}\label{mat1}
For every $i \in \{0,\ldots,n-1\}$, we have $X_i-1 \leq X_{i+1}\leq X_i+1$.
\end{claim}
\begin{subproof}
As $M_i$ is a submatrix of $M_{i+1}$ with one row and one column less than $M_{i+1}$, we have $\rk(M_i) \leq \rk(M_{i+1})\leq \rk(M_i)+2$. By definition, $X_{i+1} = X_i +1 - \rk(M_{i+1}) + \rk(M_i)$, so $X_i-1 \leq X_{i+1}\leq X_i+1$. 
\end{subproof}

\begin{claim}\label{mat2}
For any $i \in \{1,\ldots,n-1\}$ and $r \in \{1,\ldots,i\}$, we have $\Pr\left[ X_{i+1}=r-1 \mid X_i=r\right ]=1-2^{-r}$ and $\Pr\left[ X_{i+1}\geq r \mid X_i=r\right ]= 2^{-r}$.
\end{claim}
\begin{subproof}
Suppose that $X_i=r$, so $\rk(M_i)=i-r$. First let $M_i'$ be the submatrix of $M$ restricted to the first $i+1$ rows and $i$ columns. Observe that the vector space spanned by the rows of $M_i$ is of dimension $i-r$ and hence contains $2^{i-r}$ vectors. Further, the last row of $M_i'$ is chosen uniformly at random among $2^{i}$ vectors. We hence obtain $\Pr\left[\rk(M_i')=\rk(M_i)\right]=\frac{2^{i-r}}{2^{i}}=2^{-r}$ and $\Pr\left[\rk(M_i')=\rk(M_i)+1\right]= 1 -2^{-r}$. 

If $\rk(M_i')=\rk(M_i)$, then $\rk(M_{i+1}) \leq \rk(M_i) +1$, and so $X_{i+1}\geq X_i=r$.
Now consider the case $\rk(M_i')=\rk(M_i)+1$. Due to the symmetry of the matrix, we then obtain that the last column vector of $M_{i+1}$ is not spanned by the column vectors of $M_i'$. Hence $\rk(M_{i+1})=\rk(M_i')+1=\rk(M_i)+2$ yielding $X_{i+1}=(i+1)-\rk(M_{i+1})=(i+1)-(\rk(M_{i})+2)=i-\rk(M_i)-1=X_i-1$.
\end{subproof}

\begin{claim}\label{mat3}
For any $i \in \{1,\ldots,n\}$ and $r \in \{0,\ldots,i\}$, we have $\Pr\left[X_i \geq r \right]\leq 2^{-\frac{1}{2}(r^2-4r)}$.
\end{claim}
\begin{subproof}
We proceed by induction on $i$. The statement clearly holds for $i=1$. We now suppose that the statement holds for all integers up to some $i$ and show it also holds for $i+1$. For $r \leq 4$, we have $\Pr\left[X_i \geq r\right]\leq 1 \leq 2^{-\frac{1}{2}(r^2-4r)}$. We may hence suppose from now on that $r \geq 5$.
By induction, Claims~\ref{mat1} and~\ref{mat2} we obtain

\begin{eqnarray*}
\Pr\left[X_{i+1}\geq r\right]&= &\Pr\left[X_{i+1}\geq r \mbox{~and~} X_i\geq r+1\right]+\Pr\left[X_{i+1}\geq r \mbox{~and~} X_i = r\right]\\
& & 
\hspace*{7cm} + \Pr\left[X_{i+1}\geq r \mbox{~and~} X_i =r-1\right]\\
&\leq &
\Pr\left[X_i\geq r+1\right]+ \Pr\left[X_{i+1}\geq r \mid X_i = r\right]\Pr\left[X_i = r\right] \\
& & \hspace*{5cm} +\Pr\left[X_{i+1}\geq r \mid X_i = r-1\right]\Pr\left[X_i = r-1\right]\\
&\leq &\Pr\left[X_i\geq r+1\right]+2^{-r}\Pr\left[X_i=r\right]+2^{-(r-1)}\Pr\left[X_i=r-1\right]\\
&\leq & 2^{-\frac{1}{2}((r+1)^2-4(r+1))}+2^{-(r-1)}\Pr\left[X_i\geq r-1\right]\\
& = & 2^{-\frac{1}{2}((r+1)^2-4(r+1))}+2^{-(r-1)}2^{-\frac{1}{2}((r-1)^2-4(r-1))}\\
& = & 2^{-\frac{1}{2}((r^2-4r)+(2r-3))}+2^{-\frac{1}{2}((r^2-4r)+3)}\\
&\leq & 2^{-\frac{1}{2}(r^2-4r)}(2^{-\frac{1}{2}(2r-3)}+2^{-\frac{3}{2}})\\
&\leq & 2^{-\frac{1}{2}(r^2-4r)}(2^{-\frac{7}{2}}+2^{-\frac{3}{2}})\\
& \leq & 2^{-\frac{1}{2}(r^2-4r)}.
\end{eqnarray*}
\end{subproof}
The lemma now follows by applying Claim \ref{mat3} to $i=n$.
\end{proof}

For a tournament $T$ with vertex set  $[n]$ and a permutation $\sigma$ of $[n]$,
we denote by $\inv_\sigma(T)$ the minimum number of inversions needed to
transform $T$ into a transitive tournament with acyclic ordering $\sigma$.
Observe that $\inv(T) = \min_\sigma \inv_\sigma(T)$.
We also denote by $M_\sigma(T)$ the symmetric matrix in $\mathbb{F}_2^{n \times n}$ all of whose diagonal entries are 0 and which has a $0$ in a cell $(i,j)$ if and only if
the arc between $i$ and $j$ is oriented according to $\sigma$.
Finally, we say that a matrix is a diagonal matrix if all its off-diagonal entries are 0 and we let $\mathbb{D}_2^{n \times n}$ be the set of diagonal matrices in $\mathbb{F}_2^{n \times n}$.

\begin{lemma}\label{lem:inverting_with_fixed_order}
Let $T$ be a tournament on $[n]$ and let $\sigma$ be a permutation of $[n]$. Then
\[
\inv_\sigma(T) \geq \min_{D \in \mathbb{D}_2^{n \times n}}
\rk(M_\sigma(T)+D).
\]
\end{lemma}

\begin{proof}
Suppose $\inv_\sigma(T)=k$ and that ${\cal X} = (X_\ell)_{\ell\in [k]}$ is a decycling family such that $\Inv(T;{\cal X})$ is the transitive tournament with acyclic ordering $\sigma$.
For every $\ell \in [k]$, let $M_\ell$ be the matrix with cell $(i,j)$ having value
$1$ if and only if both $i$ and $j$ are in $X_\ell$.
Then if $D$ is the diagonal $n \times n$ matrix whose value at cell $(i,i)$ 
equal $1$ if and only if $i$ is in an odd number of sets in ${\cal X}$, we have
\[
M_\sigma(T) +D = M_1 + \cdots + M_k.
\]
Note that the column $j$ of each $M_{\ell}$ is either the null vector (if $j\notin X_\ell$) or the indicator vector of $X_{\ell}$ (if $j\in X_\ell$). 
Hence each $M_\ell$ has rank at most $1$. Thus $\rk(M_\sigma(T)+D) \leq \rk(M_1) + \dots + \rk(M_k) \leq k$, and the result follows.
\end{proof}

We are now ready to proceed to the proof of Theorem \ref{theorem:random_tournament}, which we recall.
\random*

\begin{proof}
Let $n$ be a positive integer and let $r \in [n]$ to be determined later.
Let $T$ be a tournament with vertex set $[n]$ taken uniformly at random.
Observe that for any fixed permutation $\sigma$ of $[n]$, the matrix
$M_\sigma(T)$ follows the uniform law in the set of symmetric matrices
in $\mathbb{F}_2^{n \times n}$ with diagonal constant to $0$.
Hence, by Lemma~\ref{lem:random_matrix}, for any fixed diagonal $n \times n$
matrix $D$ we have
\[
\Pr[\rk(M_\sigma(T)+D) \leq n-r] \leq 2^{-\frac{1}{2}(r^2-4r)}.
\]
and so by Lemma~\ref{lem:inverting_with_fixed_order} and  the Union Bound

\begin{eqnarray*}
\Pr[\inv(T) \leq n-r] 
& \leq & \sum_\sigma \sum_{D  \in \mathbb{D}_2^{n \times n}} \Pr[\rk(M_\sigma(T)+D) \leq n-r] \\
& \leq & n! 2^n 2^{-\frac{1}{2}(r^2-4r)} \\
& \leq &  2^{n \log n + n - \frac{1}{2}(r^2/2-4r)} \\
& < & 1
\end{eqnarray*}
if $n \log n + n -\frac{1}{2}(r^2-4r) < 0$.
This last condition holds if $r > 2+\sqrt{2(n\log n + n + 2)}$.
As a consequence, for every positive integer $n$, there exists a tournament of order $n$ with inversion number at least
$n-(2+\sqrt{2(n\log n + n + 2)})$ and so $\inv(n) \geq n-2\sqrt{n \log n}$ for all sufficiently large integers $n$.
\end{proof}

\section{Inversion number and dijoin}\label{sec:dijoin}

In this section, we disprove Conjecture~\ref{conj:dijoin}.
We need the following lemmas.



\begin{lemma}\label{lem:union}
Let $k \geq 3$ be an odd integer.
There is a tournament $T_k$ such that $\inv(T_k)=k$ and there are disjoint sets $A_1, \dots , A_{k-1}$ such that
$(A_1, \dots , A_{k-1}, \bigcup_{i=1}^{k-1} A_i)$ is a decycling family of $T_k$.
\end{lemma}

\begin{proof}
Let $k \geq 3$ be an odd integer and $N=4^{k-1}+1$.
For every $i \in [k-1]$, let $A_i=\{a^1_i, \dots, a^N_i\}$
and $A'_i = \{a'^1_i, \dots, a'^N_i\}$ be sets of vertices and define an ordering $<$ on these vertices by
$a^1_1 < a'^1_1 < a^2_1 < a'^2_1< \dots < a^N_1 < a'^N_1 < a^1_2 < \dots < a^N_{k-1} < a'^N_{k-1}$. We define $T_k$ as the tournament with vertex set $V(T_k) = \bigcup_{i=1}^{k-1} 
(A_i \cup A'_i)$, in which the arc between two vertices is oriented according to $<$ (i.e., $xy$ is an arc if and only if $x<y$), except the ones of the form $a^j_i a^{j'}_{i'}$ for $i \neq i'$ 
which are oriented opposite to $<$.
In other words, $T_k$ is obtained from the transitive tournament with hamiltonian path $(a^1_1, a'^1_1,  a^2_1,  a'^2_1, \dots , a^N_1 , a'^N_1 , a^1_2 , \dots , a^N_{k-1} , a'^N_{k-1})$ by reversing all the arcs between the different $A_i$.
See Figure~\ref{fig:lem_union}.

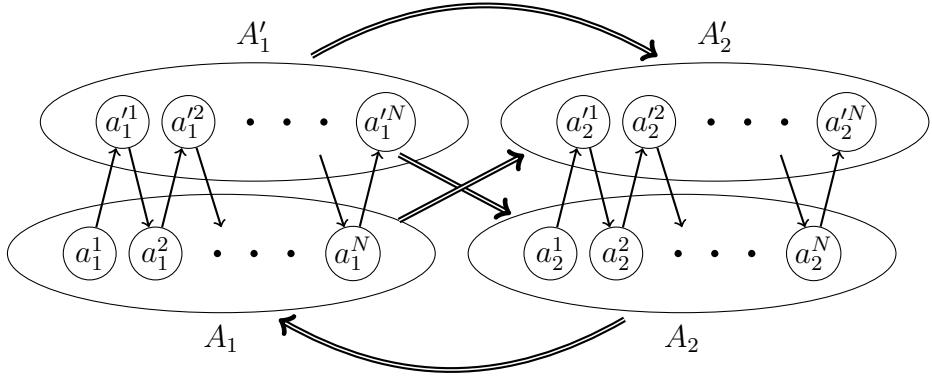
\begin{figure}[!ht]
\begin{center}
\begin{tikzpicture}[scale=1.75]
\node[vertex] (a_1_1) at (0,0) {$a_{1}^{1}$};
\node[vertex] (a'_1_1) at (0.25,1) {$a'^{1}_{1}$};
\node[vertex] (a_1_2) at (0.5,0) {$a_{1}^{2}$};
\node[vertex] (a'_1_2) at (0.75,1) {$a'^{2}_{1}$};
\node[vertex] (a_1_5) at (2.0,0) {$a_{1}^{N}$};
\node[vertex] (a'_1_5) at (2.25,1) {$a'^{N}_{1}$};
\node[vertex] (a_2_1) at (3.5,0) {$a_{2}^{1}$};
\node[vertex] (a'_2_1) at (3.75,1) {$a'^{1}_{2}$};
\node[vertex] (a_2_2) at (4.0,0) {$a_{2}^{2}$};
\node[vertex] (a'_2_2) at (4.25,1) {$a'^{2}_{2}$};
\node[vertex] (a_2_5) at (5.5,0) {$a_{2}^{N}$};
\node[vertex] (a'_2_5) at (5.75,1) {$a'^{N}_{2}$};
\draw[arc] (a_1_1) to (a'_1_1);
\draw[arc] (a'_1_1) to (a_1_2);
\draw[arc] (a_1_2) to (a'_1_2);
\draw[arc] (a'_1_2) to (1.0,0.2);
\draw[arc] (1.75,0.75) to (a_1_5);
\draw[arc] (a_1_5) to (a'_1_5);
\draw[arc] (a_2_1) to (a'_2_1);
\draw[arc] (a'_2_1) to (a_2_2);
\draw[arc] (a_2_2) to (a'_2_2);
\draw[arc] (a'_2_2) to (4.5,0.2);
\draw[arc] (5.25,0.75) to (a_2_5);
\draw[arc] (a_2_5) to (a'_2_5);
\draw (1.0, 0) ellipse (1.625 and 0.45); 
\node (label_A_1) at (1.0, -0.65) {$A_{1}$};
\draw (1.25, 1) ellipse (1.625 and 0.45); 
\node (label_A'_1) at (1.25, +1.65) {$A'_{1}$};
\draw (4.5, 0) ellipse (1.625 and 0.45); 
\node (label_A_2) at (4.5, -0.65) {$A_{2}$};
\draw (4.75, 1) ellipse (1.625 and 0.45); 
\node (label_A'_2) at (4.75, +1.65) {$A'_{2}$};
\draw[arc,double,bend left, thick] (4.0625, -0.5) to (1.4375, -0.5);
\draw[arc,double,bend left, thick] (1.6875, 1.5) to (4.3125, 1.5);
\draw[arc,double, thick] (2.35, 0.75) to (3.2, 0.3);
\draw[arc,double, thick] (2.35, 0.25) to (3.3, 0.75);

\node at (1.25, 0) {\Huge . . . };
\node at (1.5, 1) {\Huge . . . };
\node at (5, 1) {\Huge . . . };
\node at (4.75, 0) {\Huge . . . };

\end{tikzpicture}
\end{center}
\caption{The tournament $T_3$ of Lemma~\ref{lem:union}. }\label{fig:lem_union}
\end{figure}

It is straightforward to check that $(A_1, \dots , A_{k-1}, \bigcup_{i=1}^{k-1} A_i)$ is a decycling family of $T_k$.
Now we show that $\inv(T_k) \geq k$.

Suppose for a contradiction that there is a smaller decycling family ${\cal X}$ of $T_k$. Possibly adding empty sets, we may suppose that ${\cal X}=(X_1, \dots X_{k-1})$ for some $X_1,\ldots,X_{k-1}\subseteq V(T_k)$.
For every vertex $x$ in $T_k$, its {\bf indicator vector} with respect to ${\cal X}$ is the vector $\mathbf{x} = (x_1, \dots , x_{k-1}) \in \mathbb{F}_2^{k-1}$ 
such that $x_i =1$ if $x\in X_i$ and $x_i=0$ otherwise.

For every $i \in [k-1]$ and $j \in [N]$, let $\mathbf{v}^j_i$ and $\mathbf{v}'^j_i$ be the indicator vectors of $a^j_i$ and $a'^j_i$ respectively.
By the Pigeonhole Principle, for every $i \in [k-1]$, there exist
distinct indices $j,j' \in [N]$ such that $j<j'$,  $\mathbf{v}^j_i=\mathbf{v}^{j'}_i$ and
$\mathbf{v}'^j_i=\mathbf{v}'^{j'}_i$.

Set $\alpha^1_i = a^j_i$, $\alpha^2_i = a^{j'}_i$, $\alpha'^1_i = a'^j_i$, $\alpha'^2_i = a'^{j'}_i$, 
$\tilde{A}_i = \{\alpha^1_i,\alpha^2_i\}$,
$\mathbf{v}_i=\mathbf{v}^j_i=\mathbf{v}^{j'}_i$, 
$\tilde{A}'_i = \{\alpha'^1_i,\alpha'^2_i\}$ and $\mathbf{v}'_i=\mathbf{v}'^j_i=\mathbf{v}'^{j'}_i$.
We also define $\mathbf{u}_i = \mathbf{v}_i + \mathbf{v}'_i$ for every $i \in [k-1]$.
We will now extract several constraints on the vectors $\mathbf{u}_i$.

We denote by $\cdot$ the scalar product in $\mathbb{F}_2^{k-1}$. Note that $\mathbf{u} \cdot \mathbf{u} = \mathbf{1} \cdot \mathbf{u}$ for every vector 
$\mathbf{u} \in \mathbb{F}_2^{k-1}$, where $\mathbf{1} \in \mathbb{F}_2^{k-1}$ is the vector constant to $1$.

\begin{claim}\label{claim:diagonal_constraints}
    For every $i \in [k-1]$, $\mathbf{u}_i \cdot \mathbf{u}_i = 0$.
\end{claim}

\begin{subproof}
First consider the triplet $(\alpha^1_i, \alpha'^1_i,
\alpha^2_i)$.
Observe that the arcs between the pairs $\{\alpha^1_i, \alpha'^1_i\}$ and $\{\alpha'^1_i \alpha^2_i\}$ are inverted
exactly $\sum_{\ell=1}^{k-1} v_{i,\ell} v'_{i,\ell}$ times,
and the arc between the pair $\alpha^1_i \alpha^2_i$ has been inverted exactly $\sum_{\ell=1}^{k-1} v_{i,\ell}$ times.
But as $\{\alpha^1_i, \alpha'^1_i,\alpha^2_i\}$ induces originally a $TT_3$, and
must finish as a $TT_3$, we deduce that these quantities must be equal
modulo $2$.
In other words $\sum_{\ell=1}^{k-1} (v_{i,\ell} v'_{i,\ell} + v_{i,\ell}) \equiv 0 \mod 2 $, which can be rewritten $\mathbf{v}_i \cdot (\mathbf{v}'_i + \mathbf{v_i}) = 0$.

Similarly, by considering the triplet $(\alpha'^1_i, \alpha^1_i,\alpha'^2_i)$, we deduce the
constraint $\mathbf{v}'_i \cdot (\mathbf{v}'_i + \mathbf{v}_i) = 0$.
Then we deduce from the last two equalities
\begin{equation*}
    \mathbf{u}_i \cdot \mathbf{u}_i = 0, \mbox{ for all } i \in [k-1].
\end{equation*}
\end{subproof}

\begin{claim}\label{claim:outside_diagonal_constraints}
    For every distinct $i_1,i_2 \in [k-1]$, $\mathbf{u}_{i_1} \cdot \mathbf{u}_{i_2} = 1$.
\end{claim}

\begin{subproof}
Consider the triplets in $\tilde{A}_{i_1} \times \tilde{A}'_{i_1} \times \tilde{A}_{i_2} $ for distinct
$i_1, i_2 \in [k+1]$.
Observe that there is always an arc from $\tilde{A}_{i_1}$ to $\alpha'^1_{i_1}$, and an arc from $\alpha'^1_{i_1}$ to $\tilde{A}_{i_1}$. Notice that $\alpha'^1_{i_1} \ne a'^N_{i_1}$.
Further observe that all arcs between $\tilde{A}_{i_1}$ and $\tilde{A}_{i_2}$ are oriented in the same direction in $\Inv(D,\cal X)$. As $\Inv(D,\cal X)$ is acyclic, this yields that in $\Inv(D,\cal X)$, we have either $(\tilde{A}_{i_1} \cup \tilde{A}'_{i_1}) \Ra \tilde{A}_{i_2}$
or $\tilde{A}_{i_2} \Ra (\tilde{A}_{i_1} \cup \tilde{A}'_{i_1})$. Therefore,
the number of times the arcs from $\tilde{A}'_{i_1}$ to $\tilde{A}_{i_2}$ are inverted 
plus the number of times the arcs from $\tilde{A}_{i_2}$ to $\tilde{A}_{i_1}$ are inverted must be odd.
Hence $\mathbf{v}'_{i_1} \cdot \mathbf{v}_{i_2} + \mathbf{v}_{i_2} \cdot \mathbf{v}_{i_1} = 1$.
Similarly, by considering the triplets in $\tilde{A}_{i_1} \times \tilde{A}'_{i_1} \times \tilde{A}'_{i_2}$, 
we get $\mathbf{v}'_{i_1} \cdot \mathbf{v}'_{i_2} + \mathbf{v}'_{i_2} \cdot \mathbf{v}_{i_1} = 0$.
Summing these two equalities we obtain
\begin{equation*}
    \mathbf{u}_{i_1} \cdot \mathbf{u}_{i_2} = 1 ~~\mbox{~for all~} i_1,i_2 \in [k-1] \mbox{~such that~} i_1 \neq i_2.
\end{equation*}
\end{subproof}

Now consider the matrix $U \in \mathbb{F}_2^{(k-1) \times (k-1)}$ with
column number $i$ being $\mathbf{u}_i$.
Claims~\eqref{claim:diagonal_constraints} and~\eqref{claim:outside_diagonal_constraints} yield
\begin{equation}\label{eq:matrix_equation}
    U^{\top} \cdot U = 
    \left(\begin{array}{ccc}
      0 & & (1)  \\
       & \ddots & \\
      (1) & & 0 
    \end{array}\right)
    = M
\end{equation}
where $M$ is the $(k-1) \times (k-1)$ matrix 
with a $1$ in every cell, except
in the diagonal which is constant to $0$.
As $k$ is odd, we have $M^2 = I_{k-1}$.
Thus $M$ has rank $k-1$, and every solution $U$ to~\eqref{eq:matrix_equation} is invertible.
We will construct a solution $U'$ with a line equal to $0$, which is a
contradiction.
To do so, observe that in any solution $U$ to~\eqref{eq:matrix_equation}, for every column $i_1$, adding $\mathbf{1}$, the vector constant to $1$, to this column yields a new solution.
Indeed, if $\mathbf{u}_{i_1} \cdot \mathbf{u}_{i_1} = \mathbf{1} \cdot \mathbf{u}_{i_1} = 0$, $\mathbf{u}_{i_2} \cdot \mathbf{u}_{i_2} = \mathbf{1} \cdot \mathbf{u}_{i_2} = 0$ and $\mathbf{u}_{i_1} \cdot \mathbf{u}_{i_2} = 1$, 
then $(\mathbf{u}_{i_1} + \mathbf{1}) \cdot (\mathbf{u}_{i_1} + \mathbf{1}) = 0 + 0 + 0 + \mathbf{1} \cdot \mathbf{1} = 0$ (because $k-1$ is even)
and $(\mathbf{u}_{i_1} + \mathbf{1}) \cdot \mathbf{u}_{i_2} = 1 + 0$. 
Now consider $U'$ obtained from $U$ by adding the vector $\mathbf{1}$ to every column starting by a $1$.
This gives a solution to~\eqref{eq:matrix_equation} with the first row being the vector $\mathbf{0}$, which contradicts the fact that every solution to~\eqref{eq:matrix_equation}
has full rank. 

This contradiction shows that there is no decycling family of length $k-1$ for $T_k$, and so $\inv(T_k)\geq k$.
\end{proof}
We are now ready to give the proof of Theorem \ref{cor:disproof}, which we recall.

\disproof*

\begin{proof}
Let $T_k$ be as in Lemma~\ref{lem:union} with decycling family $(A_1, \dots , A_{k-1}, \bigcup_{i=1}^{k-1} A_i)$.
Let $(X_1, \dots , X_p)$ be a decycling family of $R$ with $p=\inv(R)$.
One easily checks that  $(A_1 \cup X_1, ..., A_{k-1} \cup X_1, \bigcup_{i = 1}^{k-1} A_i \cup X_1, X_2, ..., X_p)$ is a decycling family of $T_k \Ra R$.
\end{proof}

\section{Splitting into two digraphs with large inversion number}\label{sec:splitting}

In this section, we give partial positive answers to Conjecture ~\ref{pb:split} by proving Theorems \ref{cor:f-tournament} and \ref{thm:1split}.

In order to prove Theorem \ref{cor:f-tournament}, we need the following result of Belkhechine et al.~\cite{BBBP}.

\begin{lemma}\label{lem:recur}
Let $D$ be an oriented graph and let $x\in V(D)$.
Then $\inv(D) \leq \inv(D-x)+2$.
\end{lemma}

We also need a theorem of Belkhechine et al.~\cite{BBBP} on inversion-critical tournaments.
Let $k$ be a positive integer.
 A tournament $T$ is {\bf $k$-inversion-critical} if $\inv(T) =k$ and \rev{$\inv(T-x) < k$} for all $x\in V(T)$.
We denote by ${\cal IC}_k$ the set of $k$-inversion-critical tournaments.

\begin{theorem}[Belkhechine et al.~\cite{BBBP10}]\label{thm:crit-finite}
For  any positive integer $k$, the set ${\cal IC}_k$ is finite.
\end{theorem}

We are now ready to give the proof of Theorem \ref{cor:f-tournament}, which we recall.

\ftour*

\begin{proof}
For every positive integer $k$, let $p(k)$ be the maximum order of a tournament in ${\cal IC}_{k+1} \cup {\cal IC}_{k}$.
Such a number exists, by Theorem~\ref{thm:crit-finite}.

We set $f_T(k_1,k_2) = k_2 + 2p(k_1)$.

Let $T$ be a tournament with $\inv(T) \geq f_T(k_1,k_2)$.
By  Lemma~\ref{lem:recur}, $T$ has a subtournament $T_1$ in ${\cal IC}_{k_1+1} \cup {\cal IC}_{k_1}$.
 We have $\inv(T_1) \in \{k_1, k_1+1\}$ and $|T_1| \leq p(k_1)$.
Set $T_2=T-T_1$.
  By Lemma~\ref{lem:recur}, $\inv(T) \leq \inv(T_2) + 2p(k_1)$, so $\inv(T_2)\geq k_2$.
  Hence $(V(T_1), V(T_2))$ is the desired partition.
\end{proof}

We now proceed to the proof of Theorem \ref{thm:1split}. In fact, we shall prove the following result which directly implies Theorem \ref{thm:1split}.

\begin{theorem}\label{thm:f(1,k)}
Let $D$ be a non-acyclic oriented graph.
Then $D$ contains a directed cycle $C$ such that $\inv(D-C) \geq \inv(D) -  10$.
\end{theorem}
\begin{proof}
As $D$ is not acyclic, we can choose a shortest directed cycle $C=(v_1, \dots , v_p, v_1)$  in $D$.
Set $a=\lfloor p/3 \rfloor$, $b=p- 3a$.
Let $Y_1= \{v_{3i-2} \mid 1\leq i \leq a\}$, $Y_2= \{v_{3i-1} \mid 1\leq i \leq a\}$, and $Y_3= \{v_{3i} \mid 1\leq i \leq a\}$, and let
$B = V(C)\setminus (Y_1\cup Y_2\cup Y_3)$.
Note that $|B|=b \leq 2$.

First invert a minimum decycling family $\cal X$ of $D-C$.
For any $i\in \{1,2,3\}$, let $Z_i=N^+(Y_i) \setminus C$.
Observe that there is no arc from $Z_i$ to $Y_i$ for otherwise $D$ would contain a shorter directed cycle than $C$.
Therefore, after inverting $Y_i \cup Z_i$ and $Z_i$, we have all the arcs between $Y_i$ and $Z_i$ directed towards $Y_i$.
Hence after doing at most six inversions (two for each $i$) all the arcs between $V(D)\setminus V(C)$ and $Y_1\cup Y_2\cup Y_3$ are directed towards $Y_1\cup Y_2\cup Y_3$.
If $b=0$, then inverting any arc of $C$, we get a decycling family of $D$ of size $\inv(D-C) +7$. Therefore
$\inv(D-C) \geq \inv(D) -  7$.
If not, then $D\langle Y_1\cup Y_2 \cup Y_3\rangle$ is acyclic, so ${\cal X} \cup \bigcup_{i=1}^3 \{Y_i \cup Z_i, Z_i\}$ is a decycling family of $D-B$.
Hence $\inv(D-B) \leq \inv(D-C) +  6$. Thus, by Lemma~\ref{lem:recur}, $\inv(D) \leq \inv(D-B) + 2|B| \leq \inv(D-C) +  10$.
\end{proof}

\section{Further research}\label{sec:further}

\subsection{Explicit construction of tournaments with large inversion number}\label{subsec:explicit}

The lower bound $n - 2\sqrt{n\log n}$ on $\inv(n)$ is given by a random tournament. It would be interesting to have explicit constructions of tournaments with large inversion number. Theorem~\ref{thm:TT[C3]} provides tournaments of order $n$ with inversion number $\lfloor n/3\rfloor$. 
Explicit constructions of tournaments of order $n$ with inversion number larger than $\lfloor n/3\rfloor$ are not known.
An option would be to prove Conjecture~\ref{conj:Q}
Another one is to prove Conjecture~\ref{conj:triangle}, which we recall.  

\conjtriangle*

With the aid of a computer, we verified this conjecture for $k\leq 5$.
A natural way to prove it would be to prove the following extension:
{\it For any digraph $D$, $\inv(\triangle(D)) = \inv(D)+1$.}
Unfortunately, this does not hold, even for tournaments.
Consider for example the rotative $R_5$ tournament on five vertices $v_1, \dots , v_5$ in which $v_i\ra v_j$ if and only if $j-i \mod 5 \in \{1,2\}$. Let $uv$ be the arc of $\triangle(R_5) - R_5$.
Setting $X=\{v_2, v_4\}$ and $\overline{X} = \{v_1, v_3, v_5\}$, one easily checks that $\inv(R_5)=2$ and $(X, \overline{X})$  is a decycling family of $R_5$.
Thus $(X\cup\{v\}, \overline{X}\cup\{v\})$ is a decycling family of $\triangle(R_5)$.
Hence $\inv(\triangle(R_5)) = \inv(R_5) =2$.

However the statement holds for oriented graphs with inversion number at most $1$.
\begin{proposition}
If $\inv(D)\leq 1$, then $\inv(\triangle(D)) = \inv(D)+1$.
\end{proposition}
\begin{proof}
The result holds trivially when $\inv(D) =0$.

Assume now $\inv(D) =1$. Clearly, $\inv(\triangle(D)) \leq \inv(D)+1$. Suppose for a contradiction that there is a decycling set $X$ of $\triangle(D)$.
Then $X\cap V(D)$ is a non-empty strict subset of $V(D)$ for otherwise inverting $X$ would either entirely reverse $D$ or leave it unchanged. In both cases, a directed cycle would remain, a contradiction. So there is a vertex $x\in X\cap V(D)$ and a vertex $y \in V(D)\setminus X$.
Let $uv$ be the arc of $\triangle(D) - D$. By definition in $\triangle(D)$, for any $x\in V(D)$, $xu$ and $vx$ are arcs in $\triangle(D)$.
If $X\cap \{u,v\} =\emptyset$, then $(u,v, x, u)$ is a directed $3$-cycle in $\Inv(D;X)$.
If $X\cap \{u,v\} =\{u,v\}$, then $(v,u,x,v)$ is a directed $3$-cycle in $\Inv(D;X)$.
If $|X\cap \{u,v\}| =1$, then $(u,v, y, u)$ is a directed $3$-cycle in $\Inv(D;X)$.
In all cases, we get a contradiction, so $\inv(\triangle(D)) =2 = \inv(D)+1$.
\end{proof}

\begin{theorem}\label{thm:triangle-linear}
$\inv(\triangle^{k-1}(TT_2)) \geq k/2-1$.
\end{theorem}
\begin{proof}
Let $D_k = \triangle^{k-1}(TT_2)$.
We write $V(D_k) = \{a_1,b_1,a_2,b_2, \dots, a_{k},b_{k}\}$, with
$a_ib_i$ being the arc of $D_i -D_{i-1}$.
Let $\ell = \inv(D_k)$ and let ${\cal X} = (X_1, \dots X_\ell)$ be a decycling family for $D_k$.
For every $i \in [k]$, we denote by $\mathbf{v}_i$ (resp. $\mathbf{w}_i$)
the vector in $\mathbb{F}_2^\ell$ with coordinates number $j$ being $1$ if and
only if $a_i \in X_j$ (resp. $b_i \in X_j$).
Consider the matrices $V$ with columns $\mathbf{v}_i$
and $W$ with columns $\mathbf{w}_i$.
Let $M$ be the matrix defined by $M = (m_{ij})= V^\top W + I_k \in \mathbb{F}_2^{k\times k}$ where $I_k$ denotes the identity matrix of size $k$. 

\begin{claim}
$M$ contains no submatrices $\left(\begin{array}{cc}
    0 & 1 \\
    1 & 0
\end{array}\right)$ and
$\left(\begin{array}{cc}
    1 & 0 \\
    0 & 1
\end{array}\right)$.
\end{claim}
\begin{subproof}
Note that $m_{ij}$ represents the inversion of the arc 
$b_ja_i$ in $\Inv(D_k ; {\cal X})$: 
if $m_{ij}=1$, then $a_ib_j \in A(\Inv(D_k ; {\cal X}))$, and if
$m_{ij}=0$, then $b_ja_i \in A(\Inv(D_k ; {\cal X}))$.
Hence if there were a submatrix $\left(\begin{array}{cc}
    m_{i_1j_1} & m_{i_1j_2} \\
    m_{i_2j_1} & m_{i_2j_2}
\end{array}\right)$ equal to $\left(\begin{array}{cc}
    0 & 1 \\
    1 & 0
\end{array}\right)$ (resp. 
$\left(\begin{array}{cc}
    1 & 0 \\
    0 & 1
\end{array}\right)$), then $a_{i_2}b_{j_1}a_{i_1}b_{j_2}a_{i_2}$ (resp. $a_{i_1}b_{j_1}a_{i_2}b_{j_2}a_{i_1}$) would be a directed $4$-cycle in  $\Inv(D_k ; {\cal X})$, a contradiction. 
\end{subproof}

We will now show using this claim that $V^\top W = M+I$ has rank at least 
$\lceil \frac{k}{2}\rceil -1$, which will imply 
$\inv(D_k) = \ell \geq \lceil \frac{k}{2}\rceil -1$.
As $M$ contains no submatrices $\left(\begin{array}{cc}
    0 & 1 \\
    1 & 0
\end{array}\right)$ and
$\left(\begin{array}{cc}
    1 & 0 \\
    0 & 1
\end{array}\right)$, the columns of $M$ are ordered by inclusion. By sorting the columns (and the rows with the same order) according to this ordering,
we can now suppose that each row is of the form $1\dots 1 0 \dots 0$. Let $n_i$ be the number of $1$ in the row number $i$.
If $|\{i \in [k] \mid n_i\geq i \}| \geq k/2$, then the lines in $M+I$ with indices in $\{i \in [k] \mid n_i\geq i \}$ are linearly independent,
and so $M+I$ has rank at least $k/2-1$.
Otherwise, $|\{i \in [k] \mid n_i\geq i \}| < k/2$ and the lines with indices in
$\{i \in [k] \mid n_i < i \}$ in $M+I + (1)_{i,j}$ are linearly independent, and as $(1)_{i,j}$ has rank $1$, $M+I$ has rank at least $\frac{k}{2}-1$.
In both cases, $\rk(M+I) = \rk(V^\top W) \geq k/2-1$ and so $\inv(D_k) \geq k/2-1$.
\end{proof}

\subsection{Inversion number and dijoin}\label{subsec:disproof}

We disproved Conjecture~\ref{conj:dijoin}, but we believe that the conjectures implied by it, namely Conjecture~\ref{conj:tight} and Conjecture~\ref{conj:NP}, are true.

\medskip

Corollary~\ref{cor:disproof} shows that for any odd integer 
$k\geq 3$, there is a tournament $T_k$ with $\inv(T_k)=k$ such that $\inv(T_k \Ra R) < k+\inv(R)$ for all $R$ with  $\inv(R)\geq 1$.
We conjecture that the same statement  also holds for even integers.

 \begin{conjecture}
 For any $k\geq 3$, there is a tournament $T_k$ with $\inv(T_k)=k$ such that $\inv(T_k \Ra R) < k+\inv(R)$ for all $R$ with  $\inv(R)\geq 1$.
 \end{conjecture}
 
 The disproof of Conjecture~\ref{conj:dijoin} also raises the question to which extent approximate versions of this conjecture could be true.
 For any two non-negative integers $\ell$ and $r$, let $\minv(\ell,r) = \min\{ \inv(L\Ra R) \mid \inv(L) = \ell \mbox{ and }  \inv(R) = r \}$.
 \begin{problem}
Determine $\minv(\ell,r)$ for all $\ell$ and $r$.
 \end{problem}

We clearly have $\max\{\ell, r\} \leq \minv(\ell,r) \leq \ell +r$ and Corollary~\ref{cor:disproof} shows that the upper bound is often not tight.
Is the lower bound tight ?

A first problem is the following.

\begin{problem}
Does there exist a non-acyclic oriented graph $D$ such that $\inv(D \Ra D) = \inv(D)$ ?
\end{problem}

A different direction for further research is motivated by the following observation.
Every non-acyclic tournament $T$ contains a directed $3$-cycle. 
Therefore $TT_n[T]$ contains $TT_n[\vec{C_3}]$, and so
$\inv(TT_n[T]) \geq \inv(TT_n[\vec{C_3}]) =n$.
This raises the question of the following generalization of Theorem \ref{thm:TT[C3]}.
\begin{problem}
What is the maximum function $h$ such that $\inv(TT_n[T]) \geq h(\inv(T),n)$ for every tournament $T$ ?
\end{problem}

\section*{Acknowledgements}
This work was partially supported by the french Agence Nationale de la Recherche under contract Digraphs ANR-19-CE48-0013-01.
It started at a workshop at Maison Cl\'ement, Les Plantiers, devoted to digraphs problems ; we would like to express our thanks to Pierre Aboulker for organizing it (jointly with the second author).
Also, we wish to thank Pierre Charbit, Julien Duron, and Emeric Gioan, who worked with us during the workshop and later.


\end{document}